\documentclass[smallextended]{svjour3}

\newif\ifSPRINGER
\SPRINGERtrue

\usepackage{enumerate}
\usepackage[hang,flushmargin]{footmisc}
\usepackage{amsmath}
\usepackage{amsfonts}
\usepackage{amssymb}
\usepackage{url}
\usepackage{tikz,float}
\usepackage{pgfplots}
\usepackage{cite}

\usetikzlibrary{shapes}
\usepackage[affil-it]{authblk}

\newcommand*\circled[1]{\tikz[baseline=(char.base)]{
            \node[fill=black,shape=circle,draw,inner sep=1pt] (char) {\textcolor{white}{#1}};}}

\newcommand{\zero}{ 0}
\newcommand{\one}{ \circled{1}}

\titlerunning{Tighter Bounds on Directed Ramsey Number $R(7)$}

\authorrunning{D. Neiman, J. Mackey, and M.J.H. Heule}

\title{Tighter Bounds on Directed Ramsey Number $R(7)$}

\author{David Neiman, John Mackey, and Marijn Heule}

\institute{
 Carnegie Mellon University, Pittsburgh, PA, United States\\
\email{\{dneiman,jmackey,mheule\}@andrew.cmu.edu}
}

\begin{document}

\maketitle

\section*{Abstract}
Tournaments are orientations of the complete graph. The directed Ramsey number $R(k)$ is the minimum number of vertices a tournament must have to be guaranteed to contain a transitive subtournament of size $k$, which we denote by $\mathit{TT}_k$. We include a computer-assisted proof of a conjecture by Sanchez-Flores \cite{SanchezFlores54} that all $\mathit{TT}_6$-free tournaments on 24 and 25 vertices are subtournaments of $\mathit{ST}_{27}$, the unique largest $\mathit{TT}_6$-free tournament. We also classify all $\mathit{TT}_6$-free tournaments on 23 vertices. We use these results, combined with assistance from a SAT solver, to obtain the following improved bounds on $R(7)$: $34 \leq R(7) \leq 47$.


\section{Introduction and Basic Definitions}

A tournament is an orientation of a complete graph, or equivalently a directed graph $D$ with no self-loops such that, for all pairs of distinct vertices $u$ and $v$, exactly one of the edges $uv$ or $vu$ is in $D$. Intuitively, a tournament of order $n$ represents the results of a round-robin tournament between $n$ teams, where the existence of edge $uv$ means team $u$ beat team $v$ in their head-to-head match. The exclusivity of $uv$ and $vu$ means that if team $u$ beats team $v$, team $v$ doesn't beat team $u$, and the inclusion of one of those edges reflects the fact that in a round-robin tournament, each team plays each other team. The non-existence of self-loops translates to the fact that no team plays itself.

A tournament is \textit{transitive} if, for all vertices $u$, $v$, and $w$, the existence of edges $uv$ and $vw$ implies the existence of edge $uw$. The outneighborhood, or set of outneighbors, of a vertex $u$ is the set of vertices $v$ such that $u \rightarrow v$. Analogously, the inneighborhood of a vertex $u$ is the set of vertices $\{v | v \rightarrow u\}$. The definition of a tournament ensures that the inneighborhood and outneighborhood of any vertex in a tournament are disjoint, and their union includes all other vertices in the tournament.

A tournament is \textit{regular} if all vertices have the same number of outneighbors. Equivalently, this means all vertices have the same number of inneighbors, and these numbers are equal (since each edge creates one outneighbor and one inneighbor). In a round-robin tournament, this corresponds to the case where every team has the same win-loss record, which further implies that each team won exactly half of their games. Note that this can only happen if the total number of vertices is odd; for a team to win exactly half their games, there must be an even number of other teams, thus the total number of teams is odd.

A tournament is \textit{doubly regular} if it is regular and any pair $(a, b)$ of vertices has the same number of common outneighbors. These conditions imply that for any vertex $u$, its in- and out-neighborhoods are equal-sized regular tournaments. Without loss of generality, assume those neighborhoods are of size $2k + 1$, then the total number of vertices is $2(2k+1) + 1 = 4k + 3$. Thus, doubly regular tournaments can only exist when the number of vertices is congruent to 3 mod 4.

A good reference for more definitions and basic properties of tournaments is \cite{Moon}.

\section{Background}
The directed Ramsey number $R(k)$ is the smallest integer $n$ such that all tournaments on $n$ vertices contain a transitive subtournament of size $k$.

Directed Ramsey numbers were first introduced by Erd\H{o}s and Moser~\cite{ErdosMoser}. In particular, they show that $k \leq 2 \log_2(R(k)) + 1$ and also note that $k \geq \log_2(R(k)) + 1$. In particular, this means $R(k)$ grows roughly exponentially with respect to $k$, with multiplier somewhere between $\sqrt{2}$ and 2. The precise limit of that multiplier is not known, but there are known bounds on small directed Ramsey numbers:

\begin{itemize}
\item $R(2) = 2$
\item $R(3) = 4$
\item $R(4) = 8$
\item $R(5) = 14$
\item $R(6) = 28$ \cite{SanchezFlores55}
\item $32 \leq R(7) \leq 53$ \cite{SanchezFlores54,lidicky2021semidefinite}
\end{itemize}

It is also known that, for $k \leq 7$, the $\mathit{TT}_k$-free tournaments of orders $R(k)-1$ and $R(k)-2$ are unique up to isomorphism. Following the notation of Sanchez-Flores~\cite{SanchezFlores54,SanchezFlores55}, we refer to such tournaments as $\mathit{ST}_n$, where $n$ is the number of vertices in the tournament. For example, $\mathit{ST}_3$ is the unique $\mathit{TT}_3$-free tournament on 3 vertices, which is the 3-cycle.

As a side note, $\mathit{ST}_3, \mathit{ST}_7$, and $\mathit{ST}_{27}$ are the tournaments generated by quadratic residues on the appropriately-sized finite-fields (Sanchez-Flores refers to such tournaments as Galois tournaments). Galois tournaments exist for tournament sizes $n$ where $n \equiv 3$ mod $4$ and $n$ is a prime power; hence, $\mathit{ST}_{13}$ is not a Galois tournament. While Galois tournaments seem like good candidates for near-maximal Ramsey-good tournaments, it's worth noting that the Galois tournaments on 47, 43, and even 31 vertices all contain $\mathit{TT}_7$.

In this paper, we seek to improve the bounds on $R(7)$. We approach this problem, particularly the upper bound, in two steps: We first limit the set of tournaments that could possibly exist by using knowledge of smaller Ramsey numbers. We then set up the remaining cases as a Boolean satisfiability problem, which can be solved with state-of-the-art solvers.

\subsection{Boolean Satisfiability (SAT) Solvers}
The problem of finding a tournament on $n$ vertices free of $\mathit{TT}_k$s can be converted into a Boolean satisfiability problem. 
Doing so allows us to use a state-of-the-art SAT solver to show the existence of transitive tournaments. We used CaDiCaL\footnote{Commit \textsf{92d72896c49b30ad2d50c8e1061ca0681cd23e60} of\\ \url{https://github.com/arminbiere/cadical}} developed by Biere~\cite{CaDiCaL} during our experiments.

CaDiCaL, like many SAT solvers, requires that the problem be stated in conjunctive normal form (CNF). Our encodings use Boolean variables for each directed edge (where true means the edge points from the lower-index vertex to the higher-index vertex, for instance). From here, there are a few different ways to preclude the existence of a $\mathit{TT}_k$. A naive approach would be to consider all ordered combinations of $k$ vertices, negate all the edges that would create the corresponding $\mathit{TT}_k$, and require that one of those negated edges exist in the tournament; this prevents that particular $\mathit{TT}_k$. This encoding is very large and can thus only be used for small graphs or graphs for which most edges are set. 

\begin{example}
\label{ex:1}
Consider a tournament with three vertices: $u$, $v$, and $w$. For each pair of vertices $(u, v)$, we introduce a single Boolean variable that denotes the direction of the edge. For example, if the Boolean variable $uv$ is true, then there is an edge from $u$ to $v$, while if it is false, then there is an edge from $v$ to $u$. Preventing a transitive tournament of size 3 among $u$, $v$, and $w$ can be achieved by the following six clauses (all six transitive tournaments of three vertices) with $\overline{uv}$ denoting the negation of $uv$:

\begin{eqnarray*}
(uv \lor uw \lor vw) \land 
(uv \lor uw \lor \overline{vw}) \land 
(\overline{uv} \lor uw \lor vw)&\land\\
(\overline{uv} \lor \overline{uw} \lor vw) \land
(\overline{vw} \lor \overline{uw} \lor \overline{vw}) \land
(vw \lor \overline{uw} \lor \overline{vw})
\end{eqnarray*}
\end{example}

We experimented with two alternative encodings. The first one is compact and uses auxiliary variables that corresponding to the existence of a 3-cycle. For every triple of vertices $u$, $v$, and $w$, we introduce a variable $c_{\{u,v,w\}}$ that is true if and only if the edges between $u$, $v$, and $w$ form a directed cycle. From there, we need only require that for any unordered set of $k$ vertices, one of the $\binom{k}{3}$ possible 3-cycles exists. This leads to substantially smaller encoding. For large graphs with few to no edges given, this encoding is the most efficient one to solve.

\begin{example}
\label{ex:2}
Consider a graph with four vertices: $u$, $v$, $w$, and $z$. Instead of preventing a transitive tournament, we can enforce a directed cycle of size 3. We introduce a new Boolean variable $uvw$ that, if true, enforces a directed cycle between the $u$, $v$, and $w$ (in some order). The clauses that encode this constraint are:
$$
(\overline{uvw} \lor uv \lor uw) \land (\overline{uvw} \lor \overline{uv} \lor vw) \land (\overline{uvw} \lor \overline{uw} \lor \overline{vw})
\footnote{The first clause says that either there's not a cycle, or that some vertex leaves $u$ (meaning $u$ isn't the last vertex in a transitive tournament). The second and third clauses are analogous for $v$ and $w$. If $u$, $v$, and $w$ make up a transitive tournament, one of them is the last vertex; if not, they form a cycle.}$$
We can do something similar for the triples $uvz$, $uwz$, and $vwz$. In order to prevent a transitive tournament of size 4, we need to add the following additional clause: $(uvw \lor uvz \lor uwz \lor vwz)$. 
\end{example}

The final encoding that we experimented with is similar to the first one, but reduces the encoding by self-subsuming resolution~\cite{BVE}. Self-subsuming resolution works as follows: Given two clauses $C \lor x$ and $D \lor \overline x$ such that $C \subseteq D$, we can reduce $D \lor \overline x$ to $D$. This is a variant of the classical resolution rule, which would derive $C \lor D$. Since $C \subseteq D$, $C \lor D$ is logically equivalent to $D$. Self-subsuming resolution is surprisingly effective on directed Ramsey formulas and is able to reduce the first encoding by roughly 70\% on most instances that we used in our experiments. 

\begin{example}
\label{ex:3}
Continuing Example~\ref{ex:1}: The six clauses that prevent a transitive tournament of size 3 can be reduced to three clauses using self-subsuming resolution. The first two clauses, $(uv \lor uw \lor vw)$ and $(uv \lor uw \lor \overline{vw})$, can be reduced to $(uv \lor uw)$, the second two clauses  to $(\overline{uv} \lor vw)$, and the last two clauses to $(\overline{uv} \lor \overline{vw})$. These clauses can alternatively be thought of as requiring that at least one edge leaves $u$, $v$, and $w$ respectively. This is a reduction of 50\% in the number of clauses and a reduction of 67\% in the number of literals. The reduction in the number of clauses increases for larger transitive tournaments. 
\end{example}

Even after self-subsuming resolution and removing of subsumed clauses, this encoding is still larger compared to the cycle-based encoding. However, the final encoding was superior in most experiments. The main reason for the improved performance is likely that the final encoding is arc-consistent~\cite{arc-consistent}: any implication that can be made on a high-level representation, can be made by unit propagation (the main reasoning in a SAT solver). Arc-consistency is very important in SAT solvers. In Example~\ref{ex:3}, if we know that literal $uv$ is false, then unit propagation in the final encoding will deduce that variable $uw$ is true and $vw$ is false. Note that there is no unit propagation on both other encodings. 

\section{Cataloging $\mathit{TT}_6$-Free Tournaments on 24 and 25 Vertices}
\subsection{Overview}
Sanchez-Flores conjectured~\cite{SanchezFlores54} that the only $\mathit{TT}_6$-free tournaments on 24 and 25 vertices are those that appear as subtournaments of $\mathit{ST}_{27}$. He notes in his paper that this conjecture is false for 23 vertices; a counterexample can be formed by taking the circulant tournament of order 23 induced by the set of quadratic residues mod 23. These quadratic residues are the set $$\mathit{QR}_{23} = \{1, 2, 3, 4, 6, 8, 9, 12, 13, 16, 18\}$$ 

To create the circulant tournament mentioned, label the vertex set with $\{1,\dots,23\}$, and have $uv$ be an edge if $v - u \in \mathit{QR}_{23}$ mod 23. Since 23 is congruent to 3 mod 4, quadratic reciprocity tells us that $-1$ is not a quadratic residue. On the other hand, 23 is prime, so Euler's criterion tells us that there are exactly $\frac{n-1}{2}$ non-zero quadratic residues. This means for all $x \in \{1,\dots,23\}$, exactly one of $x$ and $-x$ is a quadratic residue, which is exactly the property we need for this procedure to give us a tournament. (We'll catalog 23-vertex $\mathit{TT}_6$-free tournaments more in a later section.)

In addition to being interesting on its own, proving Sanchez-Flores's conjecture would make it easier to limit the search space of $\mathit{TT}_7$-free tournaments of sizes near 50, potentially allowing us to further strengthen the upper bound on directed $R(7)$ using the techniques from the last section. We will proceed to do exactly that.

Since $\mathit{ST}_{27}$ is edge-transitive (meaning for any edges $e_1$ and $e_2$, there's an automorphism that maps $e_1$ to $e_2$), its 25-vertex subtournaments all fall into a single isomorphism class; such tournaments can be obtained by deleting any two vertices. The 24-vertex subtournaments are more interesting. One can easily generate two isomorphism classes of such subtournaments by deleting any two vertices $u$ and $v$, and then a third vertex that either formed a cycle or a $\mathit{TT}_3$ with $u$ and $v$. But we found (by deleting two fixed vertices, then looping over all remaining vertices, deleting them each in turn, and cataloging the results), that there are actually 5 isomorphism classes of 24-vertex subtournaments of $\mathit{ST}_{27}$. One of these classes is the single tournament created by deleting the unique third vertex $w$ that splits the tournament into 8 disjoint 3-cycles. (Sanchez-Flores proves existence and uniqueness of such a vertex, which is part of his proof that$\mathit{ST}_{27}$ is unique up to isomorphism~\cite{SanchezFlores55}.) The other isomorphism classes correspond to all four combinations of whether the deleted vertex $x$ formed a cycle with $u$ and $v$ and whether the deleted vertex is an outneighbor of $w$.

\subsection{Limiting the Search Space} \label{Blocks}
We now generate all $\mathit{TT}_6$-free 24-vertex tournaments. From there, 25-vertex tournaments are easy, since any $\mathit{TT}_6$-free 25-vertex tournament extends some $\mathit{TT}_6$-free 24-vertex tournament.

Let $T$ be a $\mathit{TT}_6$-free 24-vertex tournament. Without loss of generality, we choose two arbitrary vertices $u$ and $v$ from $T$ such that $u \rightarrow v$, and partition the remaining vertices into the following sets, which we'll call blocks:

\begin{itemize}
\item Let $A(u, v)$ be the set of vertices $w$ with $u \rightarrow w$ and $v \rightarrow w$.
\item Let $B(u, v)$ be the set of $w$ with $u \rightarrow w \rightarrow v$
\item Let $C(u, v)$ be the set of $w$ with $w \rightarrow u$ and $w \rightarrow v$
\item Let $D(u, v)$ be the set of $w$ with $v \rightarrow w \rightarrow u$.
\end{itemize}

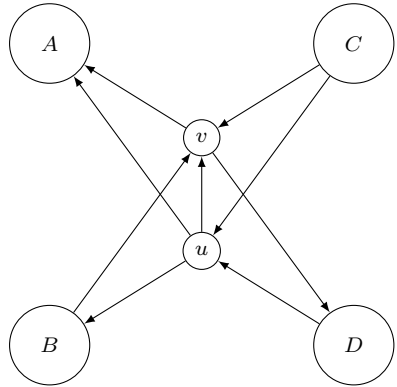
\begin{figure}
\centering
\begin{tikzpicture}

\node[circle,draw] (u) at (0,-0.75) {$u$};
\node[circle,draw] (v) at (0,0.75) {$v$};
\draw[-latex] (u) -- (v);

\node[circle,draw,minimum width=30pt] (A) at (-2,2) {$A$};
\node[circle,draw,minimum width=30pt] (B) at (-2,-2) {$B$};
\node[circle,draw,minimum width=30pt] (C) at (2,2) {$C$};
\node[circle,draw,minimum width=30pt] (D) at (2,-2) {$D$};
\draw[-latex] (u) -- (A);
\draw[-latex] (v) -- (A);
\draw[-latex] (u) -- (B);
\draw[latex-] (v) -- (B);
\draw[latex-] (u) -- (C);
\draw[latex-] (v) -- (C);
\draw[latex-] (u) -- (D);
\draw[-latex] (v) -- (D);

\end{tikzpicture}
\caption{This is a visual of the partition of $T$ into $A$, $B$, $C$, and $D$}
\end{figure}

We can now reason about the possible sizes and structures of the blocks: By adding vertices $u$ and $v$, we can extend a transitive subtournament in block $A$, $B$, or $C$ to a transitive subtournament of size two larger in $T$. This means $A$, $B$, and $C$ must be $\mathit{TT}_4$-free, which forces $|A|, |B|, |C| \leq 7$ because $R(4) = 8$.

We'll now use a cycle counting argument to limit the values of $|D|$ that we need to check: We'll find the largest possible number of 3-cycles, determine the average number of cycles on each edge, and only check $D$ block sizes up to that value.

For an arbitrary tournament $T$, the number of transitive tournaments of size 3 is just $\sum_{t \in T} \binom{d(t)}{2}$, where $t$ loops over all vertices of $T$, and $d(t)$ is the number of outneighbors of vertex $t$. This holds because each $\mathit{TT}_3$ consists of exactly one vertex of degree 2, while 3-cycles (the only other size-3 tournament) contain no vertices of degree 2. So by counting the number of ways to choose 3 vertices from $T$ with the first one dominating the other two, we exactly count the number of $\mathit{TT}_3$'s.

Since $\binom{x}{2} = x(x-1)/2 = \frac{1}{2}(x^2 - x)$ and the sum of the $d(t)$s is constant (namely $\sum_{t \in T} d(t) = \binom{n}{2}$ for a tournament on $n$ vertices, since each edge contains exactly one vertex that points to another), the number of $\mathit{TT}_3$s is minimized when $T$ is as regular as possible, and if $T$ is regular, this minimum is $n(n-1)(n-3)/8$. Since $\binom{n}{3} - \text{number of } \mathit{TT}_3\text{'s} = \text{number of 3-cycles}$, we compute that the maximum number of 3-cycles in $T$ is:
\begin{align*}
\binom{n}{3} - \frac{n(n-1)(n-3)}{8} &= \frac{n(n-1)(n-2)}{6} - \frac{n(n-1)(n-3)}{8}\\
&= \frac{4n(n-1)(n-2) - 3n(n-1)(n-3)}{24}\\
&= \frac{n(n-1) (4(n-2)-3(n-3))}{24}\\
&= \frac{n(n+1)(n-1)}{24}
\end{align*}

In this particular case, where $n = 24$ (and in particular, $n$ is even), $T$ cannot be regular, so we start by making the outdegrees as equal as possible. The average outdegree is 11.5, so we take half the vertices to have outdegree 11, and the other half to have outdegree 12. We first compute the number of $\mathit{TT}_3$'s:

\begin{align*}
\text{number of } \mathit{TT}_3 \text{'s} &= 12 \binom{12}{2} + 12 \binom{11}{2}\\
&= \frac{(12)(12)(11) + (12)(11)(10)}{2}\\
&= 1452
\end{align*}

The total number of 3-vertex subtournaments is just $$\binom{24}{3} = \frac{(24)(23)(22)}{6} = 2024$$
So the maximal number of 3-cycles is $2024 - 1452 = 572$. Each cycle contains 3 edges, so the average number of cycles on a given edge is

\begin{align*}
(572)(3)/(\text{number of edges}) &= (572)(3)/\binom{24}{2}\\
&= \frac{(572)(3)(2)}{(24)(23)}\\
&= \frac{572}{92}
\end{align*}

which is between 6 and 7. Since this average is over edges, and selecting an edge is equivalent to selecting the pair of vertices $u \rightarrow v$, we can select this edge such that the number of cycles is at most the average. Since the number of cycles is exactly the number of vertices in the $D$ block, we need only consider decompositions where $|D| \leq 6$.

We can also improve our bounds on $|A|$, $|B|$, and $|C|$ by considering the in and out neighborhoods of vertices $u$ and $v$, noting that we can extend a transitive subtournament in one of those neighborhoods by one vertex to get a transitive subtournament of $T$. In particular, this means these neighborhoods have size at most 13. It is easy to check that if a given neighborhood is of size at least 12, it cannot contain $\mathit{ST}_7$ as a subtournament (since $\mathit{ST}_{12}$ and $\mathit{ST}_{13}$ do not contain $\mathit{ST}_7$).

We now rule out the case where $A$, $B$, or $C$ has size 7. Assume for the sake of contradiction that $|A| = 7$. $\{A \cup B \cup v\}$ forms the outneighborhood of $u$, and contains an $\mathit{ST}_7$, thus has size at most 11. This means $|B| \leq 3$. Similarly, $\{A \cup D\}$ is the outneighborhood of $v$, so $|D| \leq 4$. Since $|C| \leq 7$, we have that $|A| + |B| + |C| + |D| \leq 7+3+7+4 = 21$, which is impossible because the full tournament (including $u$ and $v$) has 24 vertices, which implies $|A| + |B| + |C| + |D| = 22$. A similar argument (looking at inneighborhoods of $u$ and $v$) rules out $|C| = 7$. Finally, if $|B| = 7$, we use the outneighborhood of $u$ and the inneighborhood of $v$ to show that $|A|, |C| \leq 3$ and deduce that $|A| + |B| + |C| + |D| \leq 3+7+3+6 = 19$, which is still too small.

We now have 4 integers, each at most 6, with sum 22. The only ways for this to happen is if one of the numbers is 4 and the other three are 6, or if two of the numbers are 5 and the other two are 6. Moreover, $A$, $B$, and $C$ are $\mathit{TT}_4$-free, and $D$ is $\mathit{TT}_5$-free.

One final trick to reduce the number of cases is to consider what happens when we reverse all the edges of $T$. We can relabel $u$ as $v$ and $v$ as $u$, and then consider how the labels of $A$, $B$, $C$, and $D$ change. Vertices in $D$ still form 3-cycles with $u$ and $v$, and vertices in $B$ are still the middle vertex of $\mathit{TT}_3$s, but vertices in $A$ are now in the common inneighborhood of $u$ and $v$, while vertices in $C$ are in their common outneighborhood. Thus, by relabeling $C$ as $A$ and $A$ as $C$, we get a new decomposition of $T$ with all the edge directions flipped. Since flipping all the edges preserves whether a given set of vertices forms a transitive tournament (it simply reverses the order of such a tournament), we can without loss of generality assume $|A| \geq |C|$ assuming we also catalog the results of flipping edges of any valid tournaments obtained.

Note that flipping all the edges of $\mathit{ST}_{27}$ produces a 27-vertex tournament with no $\mathit{TT}_6$s, so that edge-flipped tournament must also be $\mathit{ST}_{27}$. This means if we only see subtournaments of $\mathit{ST}_{27}$, we don't even care about edge-flipping them, since edge-flipped subtournaments of $\mathit{ST}_{27}$ are subtournaments of the edge-flipped $\mathit{ST}_{27}$, which are still just subtournaments of $\mathit{ST}_{27}$.

Based on the above discussion, the relevant cases to test are:

\begin{table}[!h]
\centering
\begin{tabular}{l|ccccccc}
$A$ & 6 & 6 & 6 & 5 & 6 & 6 & 6\\ \hline
$B$ & 6 & 5 & 6 & 6 & 4 & 5 & 6\\ \hline
$C$ & 6 & 6 & 5 & 5 & 6 & 5 & 4\\ \hline
$D$ & 4 & 5 & 5 & 6 & 6 & 6 & 6\\
\end{tabular}
\end{table}


\subsection{Searching for Valid Tournaments}
We can now check these cases with brute-force computer search. We used Matlab for ease of debugging and visualization, and we implemented a custom constraint-propagation routine to identify impossible edge-combinations faster. This routine looks for recently-changed size-6 subtournaments with only a single unknown edge, and checks whether either edge direction would form a $\mathit{TT}_6$. If only one edge direction avoids the $\mathit{TT}_6$, that edge direction is used, and the constraint propagation repeats. If both edge directions would cause a $\mathit{TT}_6$, the current edge assignment is impossible, so we can stop early. This code is available on our GitHub repository\footnote{\url{https://github.com/neimandavid/Directed-Ramsey}}.

It's worth noting that using CaDiCaL for this problem would require rerunning the search every time a solution is found, adding blocking clauses each time to prevent finding known solutions solution (or isomorphic copies of them), until the SAT solver ultimately determines that no further solutions exist. That said, our approach isn't particularly fast either; it's possible there are more efficient ways to solve this problem.

The case when $|D| = 4$ is easy. $|A| = |B| = |C| = 6$, so $\{A \cup B\}$ and $\{B \cup C\}$ must be $\mathit{ST}_{12}$. Given that, we can (in a few hours on a laptop) determine via brute-force search that there is, up to isomorphism, only one $\mathit{TT}_6$-free way to form $T \backslash D$, and that by adding $D$ we can form any of the 5 subtournaments of $\mathit{ST}_{27}$ and nothing else.

The other cases involve a longer search, with more free variables and fewer constraints. To speed up the process, for each set of search parameters we first find all ways to connect $A$ and $B$ and all ways to connect $B$ and $C$. We then combine all compatible $AB$ and $BC$ blocks to get all possible $ABC$ blocks. Finally, we add in $D$ at the end, and find all ways to fill in all edges between $D$ and each possible $ABC$ block. We include $u$ and $v$ at all points in this process to further limit the search space. Unfortunately, $u$ and $v$ don't limit $D$ that much, which is why we add $D$ at the end and hope there aren't too many cases left to consider. We only check for isomorphisms at the end of the search, since an isomorphism between subtournaments may involve vertices moving between the $A, B$, or $C$ blocks, and the information about which vertex is in which block needs to be preserved until the end of the search.

Unsurprisingly, cases where $|D|$ is large end up being the slowest, since there are more edges to fill in during the last step, and fewer constraints to do it with. Cases where $|D| = 5$ finish in about 12 hours on a laptop. We estimated that the  case $[5, 6, 5, 6]$ would take about 10 days to do on the laptop; it finished in about 3.5 days on a cluster of computers administered by the Carnegie Mellon math department.

The remaining three cases are much slower. To get a further speedup, we parallelized the search, allowing multiple cores of the computer to each consider a subset of the $ABC$ blocks. Using a worker pool size of 16, we were able to finish the $[6, 5, 5, 6]$ case in about 2 days, and the $[6, 6, 4, 6]$ case in about 4 days. The longest case was the $[6, 4, 6, 6]$ case, which took about a week to complete. Most of this time was spent computing the $ABC$ block, since there are about 250 ways to glue an $\mathit{ST}_6$ to a 4-vertex tournament, and pairs of these combinations had to be checked. It turns out that this reduces down to only about 225 valid $ABC$ blocks, which again only lead to the five subtournaments of $\mathit{ST}_{27}$. So all 24-vertex $\mathit{TT}_6$-free tournaments are subtournaments of $\mathit{ST}_{27}$.

It is easy to verify via computer program that all of those 24-vertex $\mathit{TT}_6$-free tournaments only extend in a $\mathit{TT}_6$-free way to other subtournaments of $\mathit{ST}_{27}$.\footnote{For example, one way to check this is by calling our bruteForceEdges helper function on the appropriate partial adjacency matrices.} Since $\mathit{ST}_{27}$ is edge-transitive, it has a unique 25-vertex subtournament, obtained by deleting any two vertices. This subtournament, which we'll call $\mathit{ST}_{25}$, is thus the only $\mathit{TT}_6$-free 25-vertex tournament.

\section{Cataloguing $\mathit{TT}_6$-Free 23-Vertex Tournaments} \label{23catalog}

As noted by Sanchez-Flores~\cite{SanchezFlores54}, there exist 23-vertex $\mathit{TT}_6$-free tournaments that are not contained in $\mathit{ST}_{27}$. Using a combination of the above results, we determine that there are (up to isomorphism) exactly three such tournaments, all of which are doubly-regular.

We can use a similar argument to bound the block sizes here as we did in the 24-vertex case. A similar cycle-counting argument, which we'll use to bound the size of the $D$ block, gives:

\begin{align*}
\text{minimum number of } \mathit{TT}_3\text{'s} &= 23 \binom{11}{2} \\
&= (23)(11)(5)\\
&= (115)(11) = 1256
\end{align*}
and thus maximum number of cycles is $\binom{23}{3} - 1256 = (23)(11)(7) - (23)(11)(5) = (23)(22) = 506$. In this case, the average number of cycles on a single edge is at most

\begin{align*}
(506)(3)/\text{number of edges} &= (506)(3)/\binom{23}{2}\\
&= (23)(22)(6)/23/11 = 6
\end{align*}

and so we only need to consider cases where $|D| \leq 6$. This leaves us the following cases to check:

\begin{table}[!h]
\centering
\begin{tabular}{l|ccccccccc}
$A$ & 5 & 5 & 6 & 6 & 6 & 6 & 6 & 6 & 7\\ \hline
$B$ & 5 & 6 & 4 & 5 & 5 & 6 & 6 & 6 & 3\\ \hline
$C$ & 5 & 5 & 6 & 5 & 6 & 4 & 5 & 6 & 7\\ \hline
$D$ & 6 & 5 & 5 & 5 & 4 & 5 & 4 & 3 & 4\\
\end{tabular}
\end{table}


Since the maximum average number of cycles is exactly 6, in the $[5, 5, 5, 6]$ case we need only consider doubly-regular tournaments, since in a non-doubly-regular tournament we can choose $u \rightarrow v$ to make $|D|$ strictly less than average. McKay maintains a database of interesting combinatorial objects, including a complete list of doubly-regular 23-vertex tournaments~\cite{McKay}, so checking this case takes seconds. 3 of the 37 candidate tournaments are $\mathit{TT}_6$-free. We claim all other $\mathit{TT}_6$-free tournaments on 23 vertices are subtournaments of $\mathit{ST}_{27}$.

Most of the other cases were checked in the same way as the 24-vertex cases, using Matlab code to fill in combinations of edges, and propagating constraints to eliminate cases early. The notable exception was the $[7, 3, 7, 4]$ case. This case takes a significant amount of time to check using the Matlab code. However, in a few minutes, CaDiCaL is able to determine that no tournament with $A$, $B$, and $C$ blocks of those sizes (there's no need to even consider the $D$ block) along with the vertices $u$ and $v$, leads to an unsatisfiable result. Disappointingly, reducing any of these sizes by 1 results in a satisfiable configuration, which, combined with the slowness of the Matlab code, will make it much more time-consuming to catalog all of the 22-vertex $\mathit{TT}_6$-free tournaments.

\section{Improved Bounds on $R(7)$}
Even with our SAT representation, the space of all non-identical tournaments is too large to search with no further restrictions. However, by fixing some edges, we can make certain problems tractable and derive better bounds for $R(7)$.

\newcommand{\hs}{\hspace{0.2pt}}

\subsection{Lower Bound: $R(7) > 33$}
\begin{figure}
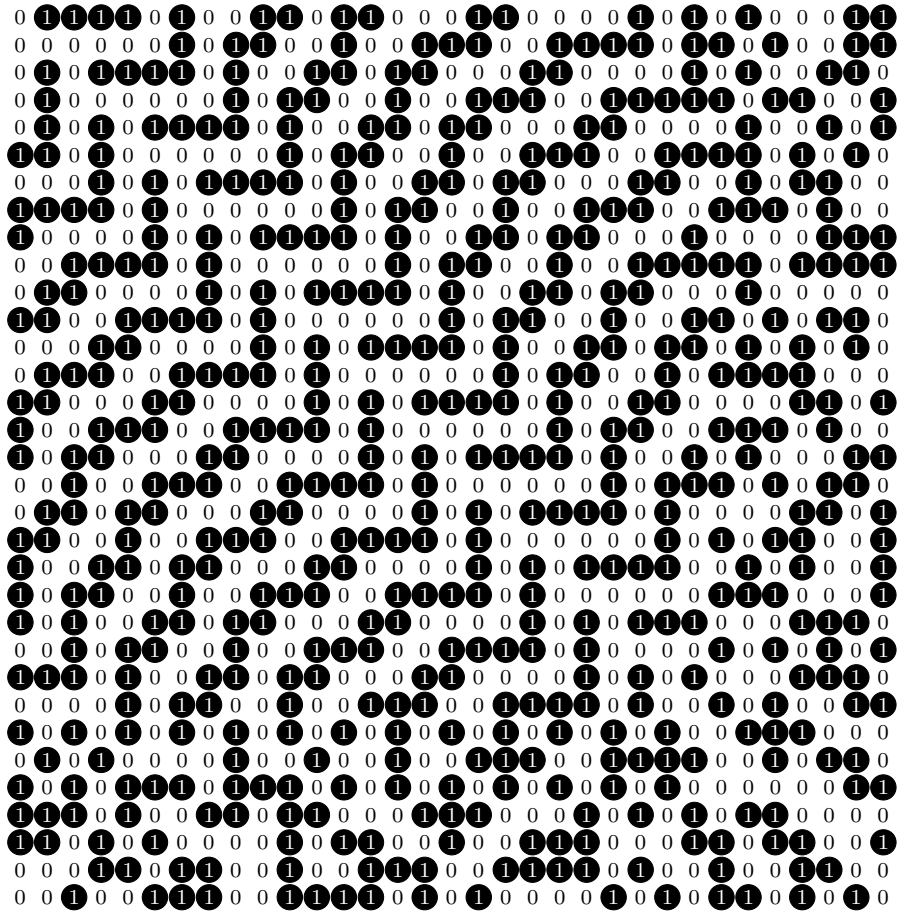

    \centering
    $$
    \begin{tabular}{@{}c@{\hs}c@{\hs}c@{\hs}c@{\hs}c@{\hs}c@{\hs}c@{\hs}c@{\hs}c@{\hs}c@{\hs}c@{\hs}c@{\hs}c@{\hs}c@{\hs}c@{\hs}c@{\hs}c@{\hs}c@{\hs}c@{\hs}c@{\hs}c@{\hs}c@{\hs}c@{\hs}c@{\hs}c@{\hs}c@{\hs}c@{\hs}c@{\hs}c@{\hs}c@{\hs}c@{\hs}c@{\hs}c@{}}
\zero& \one& \one& \one& \one& \zero& \one& \zero& \zero& \one& \one& \zero& \one& \one& \zero& \zero& \zero& \one& \one& \zero& \zero& \zero& \zero& \one& \zero& \one& \zero& \one& \zero& \zero& \zero& \one& \one\\
\zero& \zero& \zero& \zero& \zero& \zero& \one& \zero& \one& \one& \zero& \zero& \one& \zero& \zero& \one& \one& \one& \zero& \zero& \one& \one& \one& \one& \zero& \one& \one& \zero& \one& \zero& \zero& \one& \one\\
\zero& \one& \zero& \one& \one& \one& \one& \zero& \one& \zero& \zero& \one& \one& \zero& \one& \one& \zero& \zero& \zero& \one& \one& \zero& \zero& \zero& \zero& \one& \zero& \one& \zero& \zero& \one& \one& \zero\\
\zero& \one& \zero& \zero& \zero& \zero& \zero& \zero& \one& \zero& \one& \one& \zero& \zero& \one& \zero& \zero& \one& \one& \one& \zero& \zero& \one& \one& \one& \one& \one& \zero& \one& \one& \zero& \zero& \one\\
\zero& \one& \zero& \one& \zero& \one& \one& \one& \one& \zero& \one& \zero& \zero& \one& \one& \zero& \one& \one& \zero& \zero& \zero& \one& \one& \zero& \zero& \zero& \zero& \one& \zero& \zero& \one& \zero& \one\\
\one& \one& \zero& \one& \zero& \zero& \zero& \zero& \zero& \zero& \one& \zero& \one& \one& \zero& \zero& \one& \zero& \zero& \one& \one& \one& \zero& \zero& \one& \one& \one& \one& \zero& \one& \zero& \one& \zero\\
\zero& \zero& \zero& \one& \zero& \one& \zero& \one& \one& \one& \one& \zero& \one& \zero& \zero& \one& \one& \zero& \one& \one& \zero& \zero& \zero& \one& \one& \zero& \zero& \one& \zero& \one& \one& \zero& \zero\\
\one& \one& \one& \one& \zero& \one& \zero& \zero& \zero& \zero& \zero& \zero& \one& \zero& \one& \one& \zero& \zero& \one& \zero& \zero& \one& \one& \one& \zero& \zero& \one& \one& \one& \zero& \one& \zero& \zero\\
\one& \zero& \zero& \zero& \zero& \one& \zero& \one& \zero& \one& \one& \one& \one& \zero& \one& \zero& \zero& \one& \one& \zero& \one& \one& \zero& \zero& \zero& \one& \zero& \zero& \zero& \zero& \one& \one& \one\\
\zero& \zero& \one& \one& \one& \one& \zero& \one& \zero& \zero& \zero& \zero& \zero& \zero& \one& \zero& \one& \one& \zero& \zero& \one& \zero& \zero& \one& \one& \one& \one& \one& \zero& \one& \one& \one& \one\\
\zero& \one& \one& \zero& \zero& \zero& \zero& \one& \zero& \one& \zero& \one& \one& \one& \one& \zero& \one& \zero& \zero& \one& \one& \zero& \one& \one& \zero& \zero& \zero& \one& \zero& \zero& \zero& \zero& \zero\\
\one& \one& \zero& \zero& \one& \one& \one& \one& \zero& \one& \zero& \zero& \zero& \zero& \zero& \zero& \one& \zero& \one& \one& \zero& \zero& \one& \zero& \zero& \one& \one& \zero& \one& \zero& \one& \one& \zero\\
\zero& \zero& \zero& \one& \one& \zero& \zero& \zero& \zero& \one& \zero& \one& \zero& \one& \one& \one& \one& \zero& \one& \zero& \zero& \one& \one& \zero& \one& \one& \zero& \one& \zero& \one& \zero& \one& \zero\\
\zero& \one& \one& \one& \zero& \zero& \one& \one& \one& \one& \zero& \one& \zero& \zero& \zero& \zero& \zero& \zero& \one& \zero& \one& \one& \zero& \zero& \one& \zero& \one& \one& \one& \one& \zero& \zero& \zero\\
\one& \one& \zero& \zero& \zero& \one& \one& \zero& \zero& \zero& \zero& \one& \zero& \one& \zero& \one& \one& \one& \one& \zero& \one& \zero& \zero& \one& \one& \zero& \zero& \zero& \zero& \one& \one& \zero& \one\\
\one& \zero& \zero& \one& \one& \one& \zero& \zero& \one& \one& \one& \one& \zero& \one& \zero& \zero& \zero& \zero& \zero& \zero& \one& \zero& \one& \one& \zero& \zero& \one& \one& \one& \zero& \one& \zero& \zero\\
\one& \zero& \one& \one& \zero& \zero& \zero& \one& \one& \zero& \zero& \zero& \zero& \one& \zero& \one& \zero& \one& \one& \one& \one& \zero& \one& \zero& \zero& \one& \zero& \one& \zero& \zero& \zero& \one& \one\\
\zero& \zero& \one& \zero& \zero& \one& \one& \one& \zero& \zero& \one& \one& \one& \one& \zero& \one& \zero& \zero& \zero& \zero& \zero& \zero& \one& \zero& \one& \one& \one& \zero& \one& \zero& \one& \one& \zero\\
\zero& \one& \one& \zero& \one& \one& \zero& \zero& \zero& \one& \one& \zero& \zero& \zero& \zero& \one& \zero& \one& \zero& \one& \one& \one& \one& \zero& \one& \zero& \zero& \zero& \zero& \one& \one& \zero& \one\\
\one& \one& \zero& \zero& \one& \zero& \zero& \one& \one& \one& \zero& \zero& \one& \one& \one& \one& \zero& \one& \zero& \zero& \zero& \zero& \zero& \zero& \one& \zero& \one& \zero& \one& \one& \zero& \zero& \one\\
\one& \zero& \zero& \one& \one& \zero& \one& \one& \zero& \zero& \zero& \one& \one& \zero& \zero& \zero& \zero& \one& \zero& \one& \zero& \one& \one& \one& \one& \zero& \zero& \one& \zero& \one& \zero& \zero& \one\\
\one& \zero& \one& \one& \zero& \zero& \one& \zero& \zero& \one& \one& \one& \zero& \zero& \one& \one& \one& \one& \zero& \one& \zero& \zero& \zero& \zero& \zero& \zero& \one& \one& \one& \zero& \zero& \zero& \one\\
\one& \zero& \one& \zero& \zero& \one& \one& \zero& \one& \one& \zero& \zero& \zero& \one& \one& \zero& \zero& \zero& \zero& \one& \zero& \one& \zero& \one& \one& \one& \zero& \zero& \zero& \one& \one& \one& \zero\\
\zero& \zero& \one& \zero& \one& \one& \zero& \zero& \one& \zero& \zero& \one& \one& \one& \zero& \zero& \one& \one& \one& \one& \zero& \one& \zero& \zero& \zero& \zero& \one& \zero& \one& \zero& \one& \zero& \one\\
\one& \one& \one& \zero& \one& \zero& \zero& \one& \one& \zero& \one& \one& \zero& \zero& \zero& \one& \one& \zero& \zero& \zero& \zero& \one& \zero& \one& \zero& \one& \zero& \zero& \zero& \one& \one& \one& \zero\\
\zero& \zero& \zero& \zero& \one& \zero& \one& \one& \zero& \zero& \one& \zero& \zero& \one& \one& \one& \zero& \zero& \one& \one& \one& \one& \zero& \one& \zero& \zero& \one& \zero& \one& \zero& \zero& \one& \one\\
\one& \zero& \one& \zero& \one& \zero& \one& \zero& \one& \zero& \one& \zero& \one& \zero& \one& \zero& \one& \zero& \one& \zero& \one& \zero& \one& \zero& \one& \zero& \zero& \one& \one& \one& \zero& \zero& \zero\\
\zero& \one& \zero& \one& \zero& \zero& \zero& \zero& \one& \zero& \zero& \one& \zero& \zero& \one& \zero& \zero& \one& \one& \one& \zero& \zero& \one& \one& \one& \one& \zero& \zero& \one& \zero& \one& \one& \zero\\
\one& \zero& \one& \zero& \one& \one& \one& \zero& \one& \one& \one& \zero& \one& \zero& \one& \zero& \one& \zero& \one& \zero& \one& \zero& \one& \zero& \one& \zero& \zero& \zero& \zero& \zero& \zero& \one& \one\\
\one& \one& \one& \zero& \one& \zero& \zero& \one& \one& \zero& \one& \one& \zero& \zero& \zero& \one& \one& \one& \zero& \zero& \zero& \one& \zero& \one& \zero& \one& \zero& \one& \one& \zero& \zero& \zero& \zero\\
\one& \one& \zero& \one& \zero& \one& \zero& \zero& \zero& \zero& \one& \zero& \one& \one& \zero& \zero& \one& \zero& \zero& \one& \one& \one& \zero& \zero& \zero& \one& \one& \zero& \one& \one& \zero& \zero& \one\\
\zero& \zero& \zero& \one& \one& \zero& \one& \one& \zero& \zero& \one& \zero& \zero& \one& \one& \one& \zero& \zero& \one& \one& \one& \one& \zero& \one& \zero& \zero& \one& \zero& \zero& \one& \one& \zero& \zero\\
\zero& \zero& \one& \zero& \zero& \one& \one& \one& \zero& \zero& \one& \one& \one& \one& \zero& \one& \zero& \one& \zero& \zero& \zero& \zero& \one& \zero& \one& \zero& \one& \one& \zero& \one& \zero& \one& \zero\\
\end{tabular}
$$

    \caption{This is a 33-vertex tournament free of $\mathit{TT}_7$, proving that $R(7) > 33$}
    \label{33-Vertex Tournament}
\end{figure}

Using CaDiCaL, we've been able to find several $\mathit{TT}_7$-free tournaments on 33 vertices (the largest previously known was on 31 vertices), increasing the lower bound on $R(7)$. One structure that works for this is creating a vertex with inneighborhood $\mathit{ST}_{25}$ and outneighborhood 7 generic vertices (we refer to this structure as $\mathit{ST}_{25}$ extended through a pivot to 7 generic vertices). However, when we increased the outneighborhood to 8 generic vertices, CaDiCaL ran for days without reaching an answer either way. A different 33-vertex $\mathit{TT}_7$-free tournament, with possibly interesting structure, is shown in Figure~\ref{33-Vertex Tournament}. Our GitHub repository contains 84 33-vertex $\mathit{TT}_7$-free tournaments, representing 49 isomorphism classes. McKay has since used local search techniques to expand this to 5303 isomorphism classes, none of which extend to 34 vertices ~\cite{McKay}.

\subsection{Upper Bound: $R(7) \leq 47$}
First, observe the following:
\begin{itemize}
    \item $R(6) = 28$, so no indegree in a $\mathit{TT}_7$-free tournament can be larger than 27. 
    \item Reversing the directions of all edges in any tournament preserves $\mathit{TT}_7$s but swaps in- and out-neighborhoods. So if a particular size of in-neighborhood is impossible, that size of out-neighborhood is also impossible.
\end{itemize}

From here, our main tool to show the non-existence of a 47-vertex $\mathit{TT}_7$-free tournament is using CaDiCaL to rule out different in- and outdegree combinations. In the following paragraphs, we'll show the following for $\mathit{TT}_7$-free tournaments:
\begin{itemize}
    \item Indegree of at least 26 implies outdegree of at most 14 (thus no 47-vertex $\mathit{TT}_7$-free tournament can have any indegrees of 26 or 27).
    \item Indegree of 25 implies outdegree of at most 19 (thus no 47-vertex $\mathit{TT}_7$-free tournament can have any indegrees of 25)
    \item Indegree of 23 implies outdegree of at most 22 (thus no 47-vertex $\mathit{TT}_7$-free tournament can have any indegrees of 23)
    \item Each vertex of a 47-vertex $\mathit{TT}_7$-free tournament must then have indegree equal to 22 or 24. But this is impossible by a parity argument.
\end{itemize}

\subsubsection{No indegree $\geq 26$}
Using CaDiCaL, we found that $\mathit{ST}_{26}$ connected through a pivot vertex to $\mathit{TT}_5$ is unsatisfiable. Since $R(5) = 14$, any tournament containing a vertex with indegree $\geq 26$ and outdegree $\geq 14$ must contain $\mathit{TT}_7$.

\subsubsection{No indegree of 25}
Our goal here is to find a set of tournaments, such at least one must be present in any sufficiently large tournament, but such that $\mathit{ST}_{25}$ cannot extend through a pivot to any of them in a $\mathit{TT}_7$-free way.

\begin{figure}
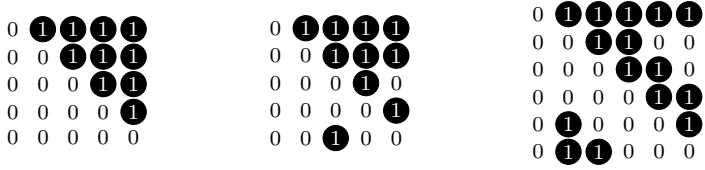

    \centering
    $
\begin{tabular}{c@{~\,}c@{\,}c@{\,}c@{\,}c}
\zero & \one & \one & \one & \one \\
\zero & \zero & \one & \one & \one \\
\zero & \zero & \zero & \one & \one \\
\zero & \zero & \zero & \zero & \one \\
\zero & \zero & \zero & \zero & \zero
\end{tabular}
$
~~~~~~~~~~
$
\begin{tabular}{c@{~\,}c@{\,}c@{\,}c@{\,}c}
\zero & \one & \one & \one & \one \\
\zero & \zero & \one & \one & \one \\
\zero & \zero & \zero & \one & \zero \\
\zero & \zero & \zero & \zero & \one \\
\zero & \zero & \one & \zero & \zero
\end{tabular}
$
~~~~~~~~~~
$\begin{tabular}{c@{~\,}c@{\,}c@{\,}c@{\,}c@{\,}c}
\zero & \one & \one & \one & \one & \one\\
\zero & \zero & \one & \one & \zero & \zero \\
\zero & \zero & \zero & \one & \one & \zero \\
\zero & \zero & \zero & \zero & \one & \one \\
\zero & \one & \zero & \zero & \zero & \one \\
\zero & \one & \one & \zero & \zero & \zero
\end{tabular}$
    \caption{The left tournament is $\mathit{TT}_5$. The middle tournament is $H_5$, the tournament consisting of a $\mathit{TT}_2$ and a 3-cycle, where the 3-cycle is in the outneighborhood of all vertices in the $\mathit{TT}_2$. The right tournament is $G_6$, the tournament consisting of a vertex with empty inneighborhood and the 5-cycle as outneighborhood. Any 10-vertex tournament contains at least one of these three tournaments as a subtournament.}
    \label{H5 G6}
\end{figure}

To start, consider the tournament on $n$ vertices consisting of a $\mathit{TT}_{n-3}$ and a 3-cycle, where the 3-cycle is in the outneighborhood of all vertices in the $\mathit{TT}_{n-3}$. Such a tournament trivially contains a $\mathit{TT}_{n-1}$, but no $\mathit{TT}_n$, and we denote this tournament $H_n$. Analogously, define $G_n$ as the tournament on $n$ vertices consisting of a $\mathit{TT}_{n-5}$ and a 5-cycle, where the 5-cycle is in the outneighborhood of all vertices in the $\mathit{TT}_{n-5}$. The adjacency matrices for $\mathit{TT}_5$, $H_5$ and $G_6$ are shown in Figure~\ref{H5 G6}.

Now consider any 10-vertex tournament $T$. By the pigeonhole principle, we can find a vertex $v$ which contains 5 vertices, say $N=\{v_1,v_2,v_3,v_4,v_5\}$, in its outneighborhood. Since the 5-cycle is the unique regular tournament on 5 vertices, if the subtournament with $N$ is the unique regular one, then $T$ has $G_6$ as a subtournament. Otherwise, again using the pigeonhole principle, $v_1$ has without loss of generality $\mathit{NN}=\{v_2,v_3,v_4\}$ among its outneighbors. If $\mathit{NN}$ is a $\mathit{TT}_3$, then $T$ contains a $\mathit{TT}_5$, and if $\mathit{NN}$ is a 3-cycle, then $T$ contains $H_5$. In any case, any tournament $T$ on 10 or more vertices must contain either $\mathit{TT}_5$, $H_5$, or $G_6$ as a subtournament.

\begin{figure}
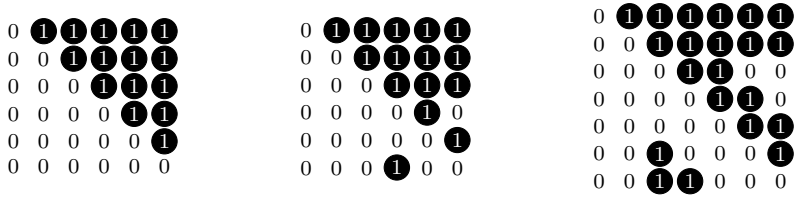

    \centering

$\begin{tabular}{c@{~\,}c@{\,}c@{\,}c@{\,}c@{\,}c}
\zero & \one & \one & \one & \one & \one \\
\zero & \zero & \one & \one & \one & \one \\
\zero & \zero & \zero & \one & \one & \one \\
\zero & \zero & \zero & \zero & \one & \one \\
\zero & \zero & \zero & \zero & \zero & \one \\
\zero & \zero & \zero & \zero & \zero & \zero
\end{tabular}$
~~~~~~~~~~
$\begin{tabular}{c@{~\,}c@{\,}c@{\,}c@{\,}c@{\,}c}
\zero & \one & \one & \one & \one & \one \\
\zero & \zero & \one & \one & \one & \one \\
\zero & \zero & \zero & \one & \one & \one \\
\zero & \zero & \zero & \zero & \one & \zero \\
\zero & \zero & \zero & \zero & \zero & \one \\
\zero & \zero & \zero & \one & \zero & \zero
\end{tabular}$
~~~~~~~~~~
$\begin{tabular}{c@{~\,}c@{\,}c@{\,}c@{\,}c@{\,}c@{\,}c}
\zero & \one & \one & \one & \one & \one & \one\\
\zero & \zero & \one & \one & \one & \one & \one\\
\zero & \zero & \zero & \one & \one & \zero & \zero \\
\zero & \zero & \zero & \zero & \one & \one & \zero \\
\zero & \zero & \zero & \zero & \zero & \one & \one \\
\zero & \zero & \one & \zero & \zero & \zero & \one \\
\zero & \zero & \one & \one & \zero & \zero & \zero
\end{tabular}$

    \caption{From left to right: $\mathit{TT}_6$, $H_6$, and $G_7$. Any 20-vertex tournament must contain one of these three as a subtournament, and $\mathit{ST}_{25}$ fails to extend through a pivot to any of these in a $\mathit{TT}_7$-free way}
    \label{Sub20}
\end{figure}

Any tournament on at least 20 vertices has a vertex with outdegree at least 10, by the pigeonhole principle. Together with the analysis of the previous paragraph, this means that any tournament on at least 20 vertices contains one of $\mathit{TT}_6$, $H_6$, or $G_7$ (see Figure~\ref{Sub20}). CaDiCaL finds that $\mathit{ST}_{25}$ does not extend in a $\mathit{TT}_7$-free way through a pivot to any of these three tournaments, which means $\mathit{ST}_{25}$ cannot be extended through a pivot in a $\mathit{TT}_7$-free way to any tournament of size at least 20.

\subsubsection{No indegree of 23}
As described in Section~\ref{23catalog}, there are two types of 23-vertex $\mathit{TT}_6$-free tournaments. There are 22 which are subtournaments of $\mathit{ST}_{25}$; there are also 3 doubly-regular tournaments. Again, our goal here is to find a large subtournament of these 23-vertex tournaments, then use CaDiCaL to show that no 23-vertex tournament extends in a $\mathit{TT}_7$-free way through a pivot to that subtournament.

We started with the subtournaments of $\mathit{ST}_{27}$. All 24-vertex $\mathit{TT}_6$-free tournaments can be formed by starting with $\mathit{ST}_{25}$ and deleting a vertex from one of five classes of vertices. Deleting one vertex from each set yields a 20-vertex tournament that, by construction, is a subtournament of any $\mathit{TT}_6$-free tournament on 24 vertices. Considering 23-vertex subtournaments of $\mathit{ST}_{27}$, we greedily deleted more vertices from this 20-vertex tournament until we found a 15-vertex tournament contained in all those 23-vertex subtournaments. We then greedily built up a 9-vertex subtournament of that 15-vertex tournament that is also contained in all three 23-vertex $\mathit{TT}_6$-free doubly-regular tournaments.

Finally, we extracted the subtournaments of that 9-vertex tournament of various sizes with the largest number of $\mathit{TT}_3$s, or equivalently the smallest number of 3-cycles; we denote by $Y_n$ the $n$-vertex subtournament obtained in this way. Given the increased difficulty of working with larger tournaments and the increased difficulty of finding $\mathit{TT}_7$-free tournaments on 32 vertices (compared to 31 vertices), we started with partial tournaments on 32 vertices, consisting of large $\mathit{TT}_6$-free tournaments extending through a pivot to some $Y_n$.

\begin{figure}
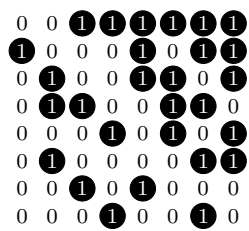

    \centering
    $$
\begin{tabular}{c@{\,}c@{\,}c@{\,}c@{\,}c@{\,}c@{\,}c@{\,}c}
\zero&\zero&\one&\one&\one&\one&\one&\one\\
\one&\zero&\zero&\zero&\one&\zero&\one&\one\\
\zero&\one&\zero&\zero&\one&\one&\zero&\one\\
\zero&\one&\one&\zero&\zero&\one&\one&\zero\\
\zero&\zero&\zero&\one&\zero&\one&\zero&\one\\
\zero&\one&\zero&\zero&\zero&\zero&\one&\one\\
\zero&\zero&\one&\zero&\one&\zero&\zero&\zero\\
\zero&\zero&\zero&\one&\zero&\zero&\one&\zero\\
\end{tabular}
$$
    \caption{This is $Y_8$, an 8-vertex tournament with few 3-cycles found in all $\mathit{TT}_6$-free 23-vertex tournaments}
    \label{Y_8}
\end{figure}

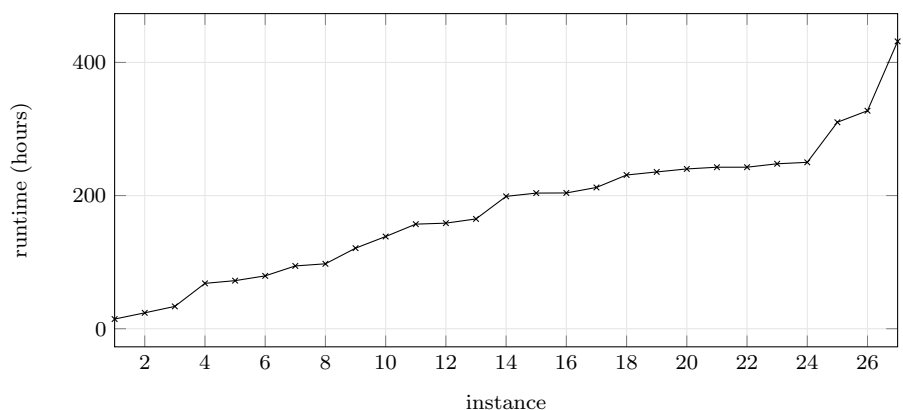
\begin{figure}[h]
\centering
   \begin{tikzpicture}
          \begin{axis}[width=\textwidth, height=6cm, mark size=1.5pt, grid=both, grid style={black!10}, xlabel=instance, ylabel=runtime (hours), xmin=1, xmax=27]

 \addplot[color=black,mark=x] coordinates {(1, 14.5195) (2, 23.8692) (3, 33.5325) (4, 68.146) (5, 72.212) (6, 79.4002) (7, 94.321) (8, 97.5169) (9, 120.995) (10, 138.503) (11, 157.089) (12, 158.58) (13, 164.923) (14, 198.812) (15, 203.735) (16, 203.937) (17, 212.188) (18, 230.885) (19, 235.513) (20, 240.008) (21, 242.591) (22, 242.607) (23, 247.749) (24, 249.936) (25, 310.129) (26, 327.438) (27, 431.605)};
 
          \end{axis}
        \end{tikzpicture}
  \caption{The runtime of CaDiCaL (in hours) of the 27 hard instances sorted by the runtime (hence the monotonic increasing curve.}
  \label{fig:runtime}
\end{figure}

We solved all 25 instances in parallel using a AWS m5d.16xlarge machine, which has 64 virtual cores and 256GB of memory.
Each instance was solved using CaDiCaL running on a single virtual core. Figure~\ref{fig:runtime} shows the runtime of these 25 instances in hours together with the two hard instances extending $ST_{25}$. Note that the average runtime was about 8 days, while one instance (23S-p-Y8) took about three weeks. The log files of the runs are available on GitHub\footnote{\url{https://github.com/neimandavid/Directed-Ramsey/tree/master/log23}}. CaDiCaL finds that none of the 25 $\mathit{TT}_6$-free 23-vertex tournaments extend through a pivot in a $\mathit{TT}_7$-free way to $Y_8$ (see Figure~\ref{Y_8}). This shows the impossibility of a $\mathit{TT}_7$-free tournament on 47 vertices possessing a vertex of indegree 23.

\subsubsection{Not all vertices can have indegree 22 or 24}
In the previous parts of this section, we eliminated the possibility that any vertex in a hypothetical 47-vertex $\mathit{TT}_7$-free tournament has indegree 23 or $\geq 25$, and the edge-flipping argument eliminates indegrees of $\leq 21$. So all the vertices in such a tournament must have indegree either 22 or 24. In particular, each indegree must be even, so the sum of the indegrees is also even. But the sum of the indegrees must equal the number of edges. Any 47-vertex tournament has $\binom{47}{2} = 1081$ edges, which is odd, a contradiction.

\section{Future Work}
It may be possible to catalog all $\mathit{TT}_6$-free tournaments with more than about 20 vertices, either using the techniques from cataloging the 23-vertex tournaments, or by using CaDiCaL or another SAT solver, and adding clauses to block anything isomorphic to each solution found. Note that while block sizes of $[7, 3, 7, 4]$ are impossible, reducing any of those sizes by 1 leads to a satisfiable problem, so we expect these cases to be significantly more difficult than the 23-vertex case. If we could create such a catalog, it could likely be used to further reduce the upper bound on $R(7)$.

On the other hand, we believe the lower bound can be increased as well. From previous patterns of a unique largest Ramsey-good tournament (that happens to be regular) and the existence of multiple 33-vertex $\mathit{TT}_7$-free tournaments, we expect that there exists at least one 35-vertex $\mathit{TT}_7$-free tournament. We've tried to get CaDiCaL to find such a tournament by starting with 17-vertex neighborhoods, but there are many  $\mathit{TT}_6$-free tournaments on 17 vertices\footnote{based on recent unpublished work by McKay, probably between about $10^{6}$ and $10^{10}$ - it's big enough that even generating a complete list of them would be hard}, and it's difficult to guess which ones are likely to be neighborhoods in a hypothetical regular 35-vertex tournament free of $\mathit{TT}_7$. Alternatively, we might increase the lower bound by incidentally finding such a 35-vertex (or larger) $\mathit{TT}_7$-free tournament while attempting to further reduce the upper bound.

\section*{Acknowledgements}
The authors would like to thank Scott Alexander (CMU undergrad) for fruitful discussions during this project, Sergei Savchenko for insightful comments about an earlier draft of this paper, and Brendan McKay for verifying and extending several of our catalogs of $\mathit{TT}_6$ and $\mathit{TT}_7$-free tournaments. 


\bibliographystyle{plain}
\bibliography{References}

\section*{Statements and Declarations}
This work was supported by NSF under grant CCF-2006363.

\section*{Data availability}
The datasets generated during the current study are available in the GitHub repository 
\url{https://github.com/neimandavid/Directed-Ramsey}.

\end{document}